\documentclass[11pt]{article}

\usepackage[utf8]{inputenc}
\usepackage{wrapfig}
\usepackage[font=small,labelfont=bf]{caption}

\usepackage[paperwidth=8.5in, paperheight=11in, margin=1.00in, footskip=20pt]{geometry}
\usepackage[font=scriptsize]{subfig}

\usepackage{graphicx,url} 
\usepackage{hyperref}
\usepackage{nopageno}
\usepackage{enumitem}
\usepackage{amsmath}
\usepackage{amsthm}
\usepackage{amssymb}
\usepackage{xcolor}
\usepackage{comment}
\usepackage{Definitions}

\title{Stochastic Rounding 2.0, with a View towards Complexity Analysis}
\author{Petros Drineas, Ilse C.F. Ipsen}
\date{}

\begin{document}

\maketitle

Stochastic Rounding (SR) is a probabilistic rounding mode that is surprisingly effective in large-scale computations and low-precision arithmetic~\cite{Croci2022}. Its random nature promotes error cancellation rather than error accumulation, resulting in slower growth of roundoff errors as the problem size increases, especially when compared to traditional deterministic rounding methods, such as rounding-to-nearest.  We advocate for SR as a foundational tool in the complexity analysis of algorithms, and suggest several research directions.

\paragraph{Stochastic Rounding 1.0.}
SR was introduced more that 70 years ago, in 1950, in a one-paragraph abstract by George Forsythe for the 52nd meeting of the \textit{American Mathematical Society}, 
and subsequently reprinted in SIAM Review \cite{Forsythe1959}. 
In the context of numerical integration
under rounding-to-nearest,
Forsythe was concerned about modeling individual round-off errors as random variables,
as first suggested by von Neumann and Goldstine in 1947.
He observed that individual roundoff errors 
at times ``had a biased distribution which caused unexpectedly
large accumulations of the [total] rounding-off error,'' and that 
``in the integration of smooth functions the ordinary rounding-off errors will
frequently not be distributed like independent random variables.''
Hence Forsythe's proposal of SR to make an individual roundoff error look like ``a true random variable.''

\paragraph{Stochastic Rounding 2.0.}
Despite its illustrious beginnings, SR has since been largely overlooked by the numerical analysis community. Yet, the hardware design landscape exhibits a strong momentum towards SR implementations on GPUs and IPUs. 

Commercial hardware, such as the Graphcore IPU, IBM floating-point adders and multipliers, AMD mixed-precision adders, and even NVIDIA's implementations of deterministic rounding, already incorporate SR in various forms,
as do the Tesla D1 and AWS Trainium chips.
SR promotes  efficient hardware integration and precision management, which are crucial for machine learning and artificial intelligence. The advent of digital neuromorphic processors like Intel Loihi and SpiNNaker2, both featuring SR, further underscores the readiness and even the \textit{necessity} for the widespread adoption of SR in hardware. The time is ripe for adoption of SR as a standard feature in GPU and IPU architectures, to optimize performance and accuracy across diverse applications.

Another reason to adopt SR is the emergence of Deep Neural Networks (DNNs) and their implications for AI. Indeed, SR is primed to be a game-changer in training neural networks with low-precision arithmetic. It tends to avoid stagnation problems typical of traditional deterministic rounding modes, thereby allowing efficient training with minimal accuracy loss. 

Prior work has demonstrated the successful training of DNNs on 16-bit and even 8-bit fixed-point arithmetic with SR, with the added benefit of a significantly reduced
energy costs. The performance of SR-based hardware compares favourably with that of higher precision formats, like binary32, but comes with a lower computational overhead. SR's application extends to dynamic fixed-point arithmetic and innovative hardware designs, such as in-memory computations and block floating-point numbers. SR has the ability to support training at various levels of arithmetic precision, and enhance convergence speeds while maintaining training accuracy. These capabilities make SR particularly attractive for AI applications where computational efficiency and resource minimization are critical.

Our piece discusses challenges and opportunities for SR in numerical linear algebra. We
 propose SR 
as a model for algorithm complexity analysis that is \textit{natively} implemented in hardware,
and discuss future research directions for alleviating the drawbacks of SR.

\paragraph{SR explained.} Suppose we want to round the number 0.7 to a single bit, either 0 or 1. Traditional deterministic rounding-to-nearest rounds 0.7 to the nearest option, which is~1. In contrast, SR rounds 0.7 to~1 with probability $0.7$ and to~0 with probability $1-0.7 = 0.3$. 
The outcome of SR rounding is a \textit{random variable} with expectation  $0.7\cdot 1 + 0.3\cdot 0 =0.7$. In statistical parlance, the rounded number is an \textit{unbiased estimator} of the exact number. This simple statement has significant, positive, implications for the behavior of SR in analyzing the numerical accuracy of arithmetic operations.

The formal definition of the stochastically
rounded version $\roundfunc(x)$ of a real
number $x$ first appeared in~\cite{parker1997monte}.
Let $\Fcal \subset \R$ be a finite set of floating point or fixed point numbers, and assume that $x$ is in the interval $[\min \Fcal, \max \Fcal]$. Identify the
two adjacent numbers in $\Fcal$ that bracket $x$,
\begin{align}\label{eqn:floorceil}
\nceil{x} = \min \{y \in \Fcal: y \geq x\} \quad \text{and}\quad \nfloor{x} = \max \{y \in \Fcal: y \leq x\},   
\end{align}
thus $\nfloor{x}\leq x\leq \nceil{x}$. 
If $\nfloor{x} = \nceil{x}$, then $\roundfunc(x) = x$, otherwise SR rounds with higher probability to the closer of the two numbers,
\begin{align*}
\roundfunc(x) = \begin{cases}\nceil{x}& 
\text{with probability}\ p(x)\equiv \frac{x - \nfloor{x}}{\nceil{x} - \nfloor{x}},\\ 
\nfloor{x} & \text{with probability}\ 1-p(x).
\end{cases}
\end{align*}
An alternative is the \textit{SR-up-or-down mode} \cite{Croci2022,vignes2004discrete}
where rounding up or down happens independently of $x$ with probability $p(x)=1/2$. 

\paragraph{Advantages of SR.} Early work~\cite{hull1966tests} gave descriptions of probabilistic models for round-off error analysis and concluded that such models are, in general, very good in theory and in practice. The recent renaissance of SR makes a strong case for randomization in the rounding process.
\begin{itemize}
\item 
SR produces roundoff errors with zero mean, thereby
encouraging cancellation of errors rather than
accumulation.
In contrast, rounding-to-nearest can
accumulate large roundoff errors whenever the individual round off errors have the same sign.

\item SR tends to avoid \textit{stagnation}.
This is a phenomenon associated
with rounding-to-nearest \cite{connolly2021stochastic,Croci2022}, where many tiny updates to a large quantity
get lost in the process of rounding.
Specifically, suppose we want compute the sum $s_0+\sum_{j=1}^k{s_j}$, where $s_0\in\Fcal$ and the magnitude of
the summands $|s_j|$ is sufficiently small compared to $|s_0|$ (smaller
than half of the distance of $s_0$ to the closest numbers in $\Fcal$), then the rounding-to-nearest operation $\fl(\cdot)$ produces
$\fl(s_0+\sum_{j=1}^k{s_j})=s_0$. This means, the 
addition of the $s_j$ does not change the sum, and the sum stagnates.

\item 
The total roundoff error from SR tends to grow more slowly with the problem size than the one from rounding-to-nearest. 
Specifically, the total
round off error from summing $n$ numbers under SR 
has, with high probability, a bound proportional to
$\sqrt{n} u$, where $u$ is the unit roundoff
\cite{hallman2023precision}. This is in contrast to 
rounding-to-nearest which has a bound proportional to~$nu$.
The proofs in \cite{hallman2023precision} rely
on measure of concentration inequalities for sums of random variables, such as Bernstein, Chernoff, and Hoeffding inequalities and corroborate the bound empirically.
\item SR increases the smallest singular value of tall-and-thin matrices, thereby improving their conditioning for the solution of least squares/regression problems. We demonstrated \cite{DBLID24} that SR implicitly regularizes such matrices in theory and in practice, therefore improving their behavior in downstream machine learning and data analysis applications. Our proofs rely
on non-asymptotic random matrix theory results which, unfortunately, lack tightness and intuition.
\end{itemize}

\paragraph{Limitations of SR.} We summarize several drawbacks of SR, as compared to rounding-to-nearest.

\begin{itemize}
    \item Lack of reproducibility: SR introduces randomness in rounding decisions, leading to non-deterministic results across different runs of the same computation. 
    \item Inability to use true randomness: Implementing SR with true randomness is impossible, therefore one has to resort to pseudo-random number generators (PRNGs). This creates an additional layer of complexity in both theoretical analysis and practical applications. 
    \item Increased computational overhead: SR requires additional computational resources to generate pseudo-random numbers and perform the rounding operation, which can slow down performance.
    \item Limited adoption in legacy systems: Existing numerical libraries and systems may not be able to support SR, limiting its adoption in legacy codes and applications.
    A detailed discussion of
    SR implementations in hardware and software can be found
    in \cite{Croci2022}.
\end{itemize}

\paragraph{Our Proposal: SR for Complexity Analysis.} The most ambitious research agenda for SR, from a numerical linear algebra and a theoretical computer science perspective, is the opportunity to establish a complexity analysis of algorithms in the presence of SR.

Existing methodologies for assessing the performance of algorithms include the following four:
\textit{Worst case analysis} upper bounds the maximum time (or space) an algorithm could possibly require, thus ensuring performance guarantees even in the most challenging of scenarios. \textit{Average case analysis} assesses the expected performance over a distribution of all possible inputs, offering a measure of efficiency for typical cases. \textit{Amortized analysis} 
 evaluates the average performance of each operation in a sequence of operations and spreads the cost of expensive operations over many cheaper ones, thus providing a balanced perspective on the overall efficiency.
\textit{Smoothed complexity}, introduced in the early 2000s \cite{Spielman2001}, aims to bridge the gap between worst and average case analyses, by evaluating algorithm performance under slight \textit{Gaussian} random perturbations of worst-case inputs, thus reflecting practical scenarios where inputs are not perfectly adversarial. Smoothed complexity has gained significant recognition for its profound impact on the analysis of algorithms and its importance is underscored by prestigious awards, including the 2008 Gödel Prize, the 2009 Fulkerson Prize, and Spielman's 2010 Nevanlinna Prize.

In contrast, we propose \textit{complexity analysis under SR}, to analyze the performance of algorithms whose operations are performed under SR. The small, random perturbations inflicted by each SR operation should, one hopes, move algorithms away from worst-case instances. In this sense, complexity analysis under SR is reminiscent of smoothed complexity but it has the major advantage of being \textit{natively} implemented in modern hardware, thus accurately reflecting algorithm performance \textit{in silico}. Therefore, SR complexity analysis has the potential to be a viable and realistic alternative to worst-case and smoothed analysis, offering a comprehensive understanding of algorithmic behavior on modern hardware. 

We can see several research directions that exploit the unique advantages of SR as a foundation for complexity analysis. First, there is a need for a theoretical framework that can rigorously define and evaluate the impact of SR on algorithm performance. A straightforward initial approach could might consist of a perturbation analysis: Analyzing the impact of SR when applied solely to the original input, with computations taking place in exact arithmetic. This approach parallels smoothed complexity analysis. A more advanced framework would extend this by evaluating the effects of SR on  the computations within the algorithm.

Second, it is crucial to establish empirical benchmarks and perform comparisons of SR against traditional deterministic rounding modes across a wide range of real-world applications. This would help in understanding the practical benefits and limitations of SR, particularly for computations on large-scale data sets in low-precision arithmetic, as is the case in ML and AI applications. These research directions could solidify the role of SR as a foundational tool in complexity analysis, akin to the role that smoothed complexity played in bridging theoretical and practical perspectives in algorithm performance assessment. 

\paragraph{Other Research Directions for SR.} We also urge the research community to focus on the directions  below:
\begin{itemize}
    \item Reproducibility in the context of SR requires pseudo-random number generators (PRNGs) with fixed seeds will be employed. While PRNGs are not truly random, they offer a well-studied, repeatable source of randomness~\cite{Barker2015}. As such, the design of PRNGs that balance performance and randomness quality optimized for hardware integration would be 
    particularly useful. SR could also be selectively applied in critical parts of the algorithm or on inputs far from rounding boundaries. More precisely, SR could be selectively employed for inputs $x$ that deviate significantly from $\nfloor{x}$ and $\nceil{x}$, as these cases result in larger perturbations and thus a more pronounced impact on the outcome.
    
    \item The development of practical non-asymptotic Random Matrix Theory tools would go a long way towards enhancing our understanding of the behavior of SR in numerical linear algebra computations. User-friendly bounds for SR would allow a better understanding of the effect of SR in standard numerical algorithms, and could lead to proofs that SR is a potent approach to reducing error propagation and increasing numerical accuracy.
    \item  A \textit{gradual} integration of SR into widely used numerical libraries and standards, combined with clear guidelines and tools for adapting legacy systems could help the incorporation of SR into existing software libraries without the need for extensive rewrites. 
\end{itemize}

\bibliographystyle{siam}
\bibliography{SR1,SR2}

\begin{thebibliography}{10}

\bibitem{Barker2015}
{\sc E.~B. Barker and J.~M. Kelsey}, {\em Recommendation for Random Number
  Generation Using Deterministic Random Bit Generators}, National Institute of
  Standards and Technology, June 2015.
\newblock NIST Special Publication 800-90A.

\bibitem{connolly2021stochastic}
{\sc M.~P. Connolly, N.~J. Higham, and T.~Mary}, {\em Stochastic rounding and
  its probabilistic backward error analysis}, SIAM J. Sci. Comput., 43 (2021),
  pp.~A566--A585.

\bibitem{Croci2022}
{\sc M.~Croci, M.~Fasi, N.~J. Higham, T.~Mary, and M.~Mikaitis}, {\em
  {Stochastic rounding: implementation, error analysis and applications}}, R.
  Soc. Open Sci., 9 (2022), p.~211631.

\bibitem{DBLID24}
{\sc G.~Dexter, C.~Boutsikas, L.~Mai, I.~C.~F. Ipsen, and P.~Drineas}, {\em
  Stochastic rounding implicitly regularizes tall-and-thin matrices},  (2024).
\newblock arXiv:2403.12278.

\bibitem{Forsythe1959}
{\sc G.~E. Forsythe}, {\em Reprint of a note on rounding-off errors}, SIAM
  Rev., 1 (1959), pp.~66--67.

\bibitem{hallman2023precision}
{\sc E.~Hallman and I.~C.~F. Ipsen}, {\em Precision-aware deterministic and
  probabilistic error bounds for floating point summation}, Numer. Math., 155
  (2023), pp.~83--119.

\bibitem{hull1966tests}
{\sc T.~E. Hull and J.~R. Swenson}, {\em Tests of probabilistic models for
  propagation of roundoff errors}, Commun. ACM, 9 (1966), pp.~108--113.

\bibitem{parker1997monte}
{\sc D.~S. Parker and D.~Langley}, {\em Monte Carlo arithmetic: exploiting
  randomness in floating-point arithmetic}, Citeseer, 1997.

\bibitem{Spielman2001}
{\sc D.~Spielman and S.-H. Teng}, {\em Smoothed analysis of algorithms}, in
  Proceedings of the thirty-third annual ACM symposium on Theory of Computing,
  ACM, 7 2001, pp.~296--305.

\bibitem{vignes2004discrete}
{\sc J.~Vignes}, {\em Discrete stochastic arithmetic for validating results of
  numerical software}, Numer. Algorithms, 37 (2004), pp.~377--390.

\end{thebibliography}
\end{document}